\documentclass[12pt,reqno]{amsart}
\usepackage{amssymb}
\usepackage{amsxtra}
\usepackage[mathscr]{eucal}

\newcommand{\Z}{{\mathbf Z}}

\newcommand{\Q}{{\mathbf Q}}
\newcommand{\N}{\mathbf{N}}

\theoremstyle{plain}

\theoremstyle{definition}

\newtheorem{Empty}{\phantom{\;\;}}
\theoremstyle{remark}

\begin{document}

\title[Counting Rationals]{Counting the Positive
Rationals: A Brief Survey}

\date{\today}

\author{David~M. Bradley}
\address{Department of Mathematics \& Statistics\\
         University of Maine\\
         5752 Neville Hall
         Orono, Maine 04469-5752\\
         U.S.A.}
\email[]{bradley@math.umaine.edu, dbradley@member.ams.org}

\subjclass{Primary: 26A03; Secondary: 11A51, 11A55, 11A63}

\keywords{Countable set, ordered pairs, enumerate, denumerable,
rational numbers, bijection, continued fraction, G\"odel
numbering, fundamental theorem of arithmetic, Egyptian fractions.}

\begin{abstract} We discuss some examples that illustrate the
countability of the positive rational numbers and related sets.
Techniques include radix representations, G\"odel numbering, the
fundamental theorem of arithmetic, continued fractions, Egyptian
fractions, and the sequence of ratios of successive hyperbinary
representation numbers.
\end{abstract}

\maketitle

%\section{Introduction}\label{sect:Intro}
Ginsberg~\cite{Gins} gave an injective mapping of the rational
numbers $\Q$ into the positive integers $\Z^+$ by interpreting the
fraction $-a/b$ as a positive integer written in base 12, with the
division slash and the minus sign as symbols for 10 and 11,
respectively. This idea was previously publicized by
Campbell~\cite{Camp}, who encountered it in the early 1970s, and
also used it to establish the countability of the ring $\Q[x]$ of
polynomials with rational coefficients. Citing Campbell and using
essentially the same arguments, Touhey~\cite{Tou} gave injections of
$\Q$ and $\Z[x]$ into $\Z^+$.  More recently, Kantrowitz~\cite{Kan}
formulated the main idea as a general principle, namely that the set
of all words of finite length that may be formed from the letters of
a countable alphabet is countable.

If we identify the set $\Q^+$ of positive rational numbers with the
subset of $\Z^+\times\Z^+$ consisting of the ordered pairs of
coprime positive integers, then any injective mapping
$\varphi:\Z^+\times\Z^+\to\Z^+$ restricts to an injection of $\Q^+$
into $\Z^+$.  Additional examples follow.  All maps $\varphi$ have
domain $\Z^+\times\Z^+$ and codomain $\Z^+$.

\begin{Empty}\label{oddpower2} An extremely simple \emph{bijection} is
given by  $\varphi(n,m)=(2n-1)2^{m-1}$. Bijectivity of $\varphi$ is
equivalent to the fact that every positive integer has a unique
representation as the product of an odd positive integer and a power
of 2.
\end{Empty}

\begin{Empty}\label{base1} Let $\varphi(n,m)=(2^n-1)2^m$.  It is easy to verify
that $\varphi$ is injective using elementary divisibility
properties of the positive integers.  However, if $(2^n-1)2^m$ is
written in binary (base 2), then $\varphi$ can be regarded as the
map which encodes the ordered pair $(n,m)$ as the string
$\{1\}^n\{0\}^m$ of $n$ consecutive ones followed by $m$
consecutive zeros. From this viewpoint, injectivity is obvious.

Grant and Priest~\cite{Gran} offered the surjection $g:\Z^+\to\Q^+$
defined by $g(x)=(k+1)/(m+1)$, where $k$ and $m$ are the number of
fours and sevens, respectively, in the decimal representation of
$x$.  This is not entirely unrelated to our example, as changing the
words ``four'' to ``one'', ``seven'' to ``zero'' and ``decimal'' to
``binary'' reveals.
\end{Empty}

\begin{Empty}\label{base2} Let $\varphi(n,m) = 2^{n+m+1}-1-2^m$. Again, it
is easy to verify directly that $\varphi$ is injective, and again
there is a nice interpretation of the map via base change.
Essentially, $\varphi$ encodes the ordered pair $(n,m)$ by writing
$n$ and $m$ in unary (base 1) using the symbol 1, separating the
two strings of ones with a zero, and then interpreting the
resulting string $\{1\}^n0\{1\}^m$ as the binary (base 2)
representation of the positive integer $2^{n+m+1}-1-2^m$ in the
usual way.
\end{Empty}

\begin{Empty}\label{base2a} A simpler injection from the viewpoint
of the algebraic formulation of the rule is
$\varphi(n,m)=2^n+2^{n+m}$, which translates to the binary string
$1\{0\}^{m-1} 1\{0\}^n$.
\end{Empty}

\begin{Empty}\label{Godel}  Let $\varphi(n,m) = 2^n3^m$.
Injectivity is a consequence of the fundamental theorem of
arithmetic.  This is a very simple example of what is sometimes
referred to as a G\"odel numbering, in which various elements are
mapped to products of prime powers.  The technique can be easily
adapted to establish countability of a wide variety of countable
sets, such as the set of all finite subsets of a countable set, a
countable union of countable sets, and so on. In~\cite{Frank} it is
used to establish countability of the ring $\Q[x]$ of polynomials
with rational coefficients, an approach that contrasts with
Campell's~\cite{Camp}.
\end{Empty}

\begin{Empty}\label{Priest} Let
$\mathscr{P}=\{2,3,5,7,11,\dots\}$ be the set of prime positive
integers and let $p:\Z^+\to\mathscr{P}$ be injective. Then
$\phi(n,m)=(p(m))^n$ defines an injection of $\Z^+\times\Z^+$ into
$\Z^+$. The related map $(p(m))^n\mapsto n/m$ is extended to a
surjection of $\Z^+$ onto $\Q^+$ in~\cite{Gran}.
\end{Empty}

The injections of examples~\ref{base1} through~\ref{Priest} grow too
quickly to be surjective.  However, if the growth rate is reduced
from exponential to quadratic, it is possible to get bijective maps
distinct from the bijection of example~\ref{oddpower2}.

\begin{Empty}\label{Cantor} Let $\varphi(n,m) = (n+m-1)(n+m-2)/2+n$. It is
an interesting exercise to prove algebraically that $\varphi$ is
bijective.  By following Cantor, a visual proof is readily obtained.
Regard $\Z^+\times\Z^+$ as an infinite-dimensional matrix with
$(n,m)$ in row $n$ and column $m$.  In light of the identity
$1+2+\cdots+(n+m-2)=\varphi(n,m)-n$, one sees that $\varphi$ lists
the entries in order starting with $(1,1)$, $(1,2)$, $(2,1)$,
$(1,3)$, $(2,2)$, $(3,1)$ and   traversing successive diagonals with
$n+m$ constant so that as each diagonal is traversed, the row index
increases as the column index decreases.
\end{Empty}

\begin{Empty}\label{Ls} Picture $\Z^+\times\Z^+$ as the lattice of
points with positive integer coordinates in the first quadrant of
the Cartesian plane.  Instead of traversing diagonals, we can
exhaust the lattice by tracing out successively larger upside-down
capital L's, starting with $(1,1)$, then $(2,1)$, $(2,2)$, $(1,2)$,
then $(3,1)$, $(3,2)$, $(3,3)$, $(2,3)$, $(1,3)$, and so on.  In
this scheme, the lattice point $(n,m)$ occurs at position
$\varphi(n,m)=(\max(n,m))^2-\max(n,m)+m-n+1$ in the sequence. Again,
it is an interesting exercise to prove algebraically that $\varphi$
is bijective.

Alternatively, one can start by listing the pairs $(n,m)$ with $n\le
m$ (the horizontal portion of the upside-down L's), and then insert
into every other position the remaining pairs obtained by switching
$n$ and $m$ when they differ.  This gives a list that begins
\begin{equation}\label{L2}
(1,1);\; (1,2), (2,1), (2,2);\; (1,3), (3,1), (2,3), (3,2), (3,3);
\end{equation}
and so on.

MacHale~\cite{Hale} gave a visual bijective correspondence between
the complete integer lattice $\Z\times\Z$ and the positive integers
using a spiral path starting at the origin.  With a little effort,
one should be able to provide an algebraic formula for  the position
of a generic lattice point in the sequence under this scheme as
well. Other more complicated ``array-based'' enumerations may be
found in~\cite{God,Har,Jon}.
\end{Empty}

\begin{Empty} An explicit bijective correspondence between $\Q^+$ and $\Z^+$ can be
obtained by exploiting the multiplicative structure of these sets
more fully. Again, let $\mathscr{P}=\{2,3,5,7,11,\dots\}$ be the set
of prime positive integers. By the fundamental theorem of
arithmetic, the map
\begin{equation}\label{Q^+toZ^*}
   \prod_{p\in\mathscr{P}} p^{\alpha_p} \mapsto
   \tilde{\alpha}=(\alpha_2,\alpha_3,\alpha_5,\alpha_7,\dots)
\end{equation}
defines a bijection from $\Q^+$ to the set $\Z^*$ of integer
sequences indexed by $\mathscr{P}$, all but finitely many terms of
which are zero.  Let $\N$ denote the set of non-negative integers,
and let $\N^*$ denote the set of all sequences of non-negative
integers indexed by $\mathscr{P}$, all but finitely many terms of
which are zero.  If $\tilde{\alpha}\mapsto\tilde{\beta}$ is any
bijection from $\Z^*$ to $\N^*$, then
\begin{equation}\label{alphabeta}
   \prod_{p\in\mathscr{P}} p^{\alpha_p} \mapsto
   \prod_{p\in\mathscr{P}} p^{\beta_p}
\end{equation}
defines a bijection from $\Q^+$ to $\Z^+$.

For an example of a bijection $\tilde{\alpha}\mapsto \tilde{\beta}$,
one could take for each $p\in\mathscr{P}$ any bijection
$g_p:\Z\to\N$ such that $g_p(0)=0$, and set $\beta_p=g_p(\alpha_p)$.
But it's probably simplest to use the same bijection $g_p=g:\Z\to\N$
for each coordinate. One possible choice for $g$ is the map
\[
   \sum_{k=0}^\infty a_k (-2)^k \mapsto \sum_{k=0}^\infty a_k 2^k,
\]
where each $a_k\in\{0,1\}$ and all but finitely many are zero.
Another possibility is to let $g(\alpha)=2\alpha$ if $\alpha\ge 0$
and $g(\alpha)=-2\alpha-1$ if $\alpha<0$.  With this choice, if for
each $p\in\mathscr{P}$, $\beta_p=g(\alpha_p)$ in~\eqref{alphabeta},
then $3/5$ maps to $3^25=45$ and the two millionth positive rational
number is $2^{-4}5^3=125/16$.  This is precisely the map presented
in~\cite{Sag}; earlier we find the inverse map given in~\cite{Frei},
but the idea of a bijection of the form~\eqref{alphabeta} goes back
earlier still---at least to~\cite{Mc}. In each instance
after~\cite{Mc} the authors, journal editors and referees all seem
to have been unaware that the same idea appeared previously.

Niven's bijective correspondence~\cite{Niv} between the rationals
and the positive integers exploits the representation for the
positive rationals as products of primes with integer exponents, as
shown here on the left hand side of~\eqref{alphabeta}.  In light of
this and the fact that he was a number theorist, it seems remarkable
that he failed to exploit the corresponding (and if anything, more
familiar) canonical representation of positive integers as products
of prime powers, shown here on the right hand side
of~\eqref{alphabeta}. Instead, he set up an intermediate
correspondence between integers written in binary notation and
sequences derived from examining consecutive blocks of ones in the
binary representation.  Thus, his bijection is unnecessarily
complicated.
\end{Empty}

\begin{Empty}
Let $p_0,p_1,p_2,\dots$ be any enumeration of the prime positive
integers (i.e.\ the map $n\mapsto p_n$ defines a bijection from
$\N$ to $\mathscr{P}$).  If we index the elements of $\Z^*$ by
$\N$ instead of $\mathscr{P}$, then the map~\eqref{Q^+toZ^*}
becomes
\[
   \prod_{n\in\N} p_n^{\alpha_n} \mapsto \tilde{\alpha} =
   (\alpha_0,\alpha_1,\alpha_2,\dots).
\]
Let $x$ be an indeterminate.  Forming the generating function of
the sequence $\tilde{\alpha}$ gives a function from $\Q^+$ to
$\Z[x]$ defined by
\[
   \prod_{n\in\N} p_n^{\alpha_n} \mapsto \sum_{n\in\N} \alpha_n
   x^n,
\]
which is not merely a bijective set map, but in fact a \emph{group
isomorphism}~\cite[II.1, p.\ 75, ex.\ 11]{Hung} between the
multiplicative free abelian group $(\Q^+,\cdot)$ with basis
$\mathscr{P}$ and the additive free abelian group $(\Z[x],+)$ with
basis equal to the non-negative integer powers of $x$.  The
abelian monoids $(\Z^+,\cdot)$ and $(\N[x],+)$ are likewise
isomorphic.  Thus, given any bijection $\sum_{n\ge 0} \alpha_n x^n
\mapsto \sum_{n\ge 0} \beta_n x^n$ from $\Z[x]$ to $\N[x]$, we get
a bijection from $\Q^+$ to $\Z^+$ via
\[
   \prod_{n\in\N} p_n^{\alpha_n} \mapsto
   \sum_{n\in\N} \alpha_n x^n \mapsto
   \sum_{n\in\N} \beta_n x^n \mapsto
   \prod_{n\in\N} p_n^{\beta_n}.
\]

\end{Empty}

\begin{Empty}  The theory of continued fractions yields additional
possibilities.  It is known~\cite{HardyWright} that every
non-negative rational number $r$ has a unique representation as a
finite continued fraction of the form
\[
   r = a_0 + \cfrac{1}{a_1+\cfrac{1}{a_2+ \cfrac{1}{{\ddots} +
   \cfrac{1}{a_{n-1}+\cfrac{1}{a_n+\cfrac{1}{1}}}}}},
\]
where $a_0,a_1,\dots,a_n$ are integers with $a_0\ge 0$ and $a_j\ge
1$ for $1\le j\le n$.  Therefore, if $b_k=\sum_{j=0}^k a_j$ %$b_k=a_0+a_1+\cdots+a_k$
for $0\le k\le n$, then
\[
   r \mapsto \sum_{k=0}^n 2^{b_k}
\]
defines a bijection from the set $\Q_{\ge 0}$ of non-negative
rational numbers to the set $\Z^+$ of positive integers. A more
complicated bijective correspondence between $\Q_{\ge 0}$ and $\Z^+$
using continued fractions is described in~\cite{Czyz}.
\end{Empty}

\begin{Empty} Of course, a bijection from $\Z^+$ to $\Q^+$
may viewed as a sequence which enumerates each positive rational
exactly once.  The recursion~\cite{Ting}
\begin{equation}\label{Ting}
   \gamma_1=1,\qquad \gamma_{2k}=1+\gamma_k,\qquad
   \gamma_{2k+1}=1/\gamma_{2k},\qquad k\in\Z^+
\end{equation}
defines one such sequence.  We briefly reproduce the motivation
for this definition here.  Start by enumerating the finite
non-empty sets of positive integers: let $M_1=\{1\}$ and for $k\ge
1$, let $M_{2k}=\{n+1: n\in M_k\}$, $M_{2k+1}=M_{2k}\cup \{1\}$.
An easy induction argument proves that for each positive integer
$k$, $\{M_j : 1\le j < 2^k\}$ coincides with the collection of
non-empty subsets of $\{1,2,\dots,k\}$.  Next, define an
enumeration of all finite sequences of positive integers in such a
way that $S_k=(a_1,a_2,\dots,a_n)$ if and only if
$M_k=\{a_1,a_1+a_2,\dots,a_1+a_2+\cdots+a_n\}$.  Finally, for
$k\in\Z^+$, let
\[
   \rho_k =\cfrac{1}{a_1+\cfrac{1}{a_2+\cfrac{1}{{\ddots}+\cfrac{1}{a_n+\cfrac{1}{1}}}}}
\]
and $\gamma_k=1/\rho_k-1$.  Then $\gamma_1=1$, and if $n>1$, then
\[
   \gamma_k =
   a_1-1+\cfrac{1}{a_2+\cfrac{1}{a_3+\cfrac{1}{{\ddots}+\cfrac{1}{a_n+\cfrac{1}{1}}}}}.
\]
Since $S_{2k}$ is just $S_k$ with $a_1$ replaced by $a_1+1$ and
$S_{2k+1}$ is just $S_k$ with an extra 1 in front, the
recurrence~\eqref{Ting} follows immediately.  The fact that every
positive rational occurs precisely once in the sequence $\gamma_1,
\gamma_2, \gamma_3,\dots$ follows from the existence and uniqueness
of the continued fraction representation.
%let $M_0=\{\}$, $M_1=\{1\}$,
%$M_2=\{2\}$, $M_3=\{1,2\}$ and in general, if $M_j$ has been defined
%for $0\le j<2^k$, then let $M_{2^k+j} = M_j \cup \{k+1\}.$
\end{Empty}

\begin{Empty} Lauwerier~\cite[p.\ 23]{Law} describes a listing of the positive
rationals between 0 and 1 arranged first by increasing denominator
and then by increasing numerator:
\[
   \frac11,\quad \frac12,\quad \frac13,\quad \frac23,\quad \frac14,\quad
   \frac34,\quad
   \frac15,\quad \frac25,\quad \frac35,\quad \frac45,\quad \frac16,\quad \frac56,\quad \dots
\]
A complete list of the positive rationals is obtained by inserting
the reciprocals:
\[
  \frac11,\quad \frac12,\quad  \frac21,\quad \frac13,\quad  \frac31,\quad \frac23,\quad
  \frac32,\quad \frac14,\quad \frac41,\quad \frac34,\quad \frac43,\quad \frac15,\quad
   \frac51,\quad\dots
\]
This list can also be obtained from~\eqref{L2} by sending the
ordered pair $(n,m)$ to $n/m$, omitting pairs in which $n$ and $m$
have a common divisor greater than 1.
\end{Empty}

\begin{Empty} An enumeration of the positive rationals with
combinatorial significance is discussed in~\cite{CalkWilf}. Let
$b_0=1$ and for $n\in\Z^+$, let $b_n$ be the number of
\emph{hyperbinary} representations of $n$.  That is, $b_n$ is the
number of ways to write $n$ as a sum of powers of 2, each power
being used at most twice. For example, $b_5=2$ because
$5=2^2+2^0=2^1+2^1+2^0.$ It's not hard to see that for all $k\ge 0$,
$b_{2k+1}=b_k$ and $b_{2k+2}=b_k+b_{k+1}$. This recursion together
with the initial condition gives an alternative way of defining the
sequence of hyperbinary representation numbers; it also provides a
convenient way to compute the first several terms. Calkin and Wilf
proved that each positive rational number occurs once and only once
in the list $b_0/b_1, b_1/b_2, b_2/b_3,\dots$ Thus, the map
$\psi:\Z^+\to\Q^+$ defined by $\psi(k)= b_{k-1}/b_k$ is bijective.
\end{Empty}

\begin{Empty}
Cohen~\cite{Coh} proves inter alia that every positive rational $r$
such that $0<r<1$ has a unique representation as a sum of unit
fractions of the form
\[
   r = \sum_{j=1}^k \prod_{i=1}^j \frac{1}{n_i},
\]
where the $n_i$ are integers and $2\le n_1\le n_2\le \cdots \le
n_k$.  If we let $\psi(r) = \sum_{j=1}^k 2^{n_j+j-3}$, then the map
$x\mapsto \psi(x/(x+1))$ defines a bijection from $\Q^+$ to $\Z^+$.
\end{Empty}

\end{document}